\documentclass[12pt,a4paper]{amsart}
\usepackage{geometry}
\usepackage{tabls}
\usepackage{url}
\usepackage[utf8]{inputenc}
\usepackage{amsfonts}
\usepackage{amssymb}
\usepackage{mathrsfs}
\usepackage{amsmath}
\usepackage{amsthm}
\usepackage{amscd}
\usepackage{fancyhdr}
\usepackage{enumerate}
\usepackage{paralist}
\usepackage{xypic}
\usepackage{graphicx}
\usepackage{cite}
\usepackage{lipsum}
\usepackage{rotating}
\usepackage{footmisc}
\usepackage[all,cmtip]{xy}

\numberwithin{equation}{section}

\theoremstyle{plain}
\newtheorem{Cor}{Corollary}[section]
\newtheorem{Def}[equation]{Definition}
\newtheorem{Thm}[equation]{Theorem}

\newtheorem{prop}[equation]{Proposition}
\newtheorem{rem}[equation]{Remark}
\newtheorem{ex}[equation]{Example}

\begin{document}

\title{Iterated integrals on products of one variable  multiple polylogarithms}

\author{Jiangtao Li}

\email{lijiangtao@amss.ac.cn}
\address{Jiangtao Li \\Hua Loo-Keng Center for Mathematics Sciences,
          Academy of Mathematics and Systems Science,
         Chinese Academy of Sciences, 
         Beijing, China}

\begin{abstract}
    In this paper, we show that the iterated integrals on products of one variable multiple polylogarithms from $0$ to $1$ are actually in the algebra of multiple zeta values if they are convergent. In the divergent case, we define the regularized iterated integrals from $0$ to $1$. By the same method, we show that the regularized iterated integrals are also in the algebra of multiple zeta values. As an application,  we  give new series representations for  multiple zeta values and calculate some interesting examples of iterated integrals.
     
\end{abstract}

\let\thefootnote\relax\footnotetext{
Project funded by  China Postdoctoral Science Foundation grant 2019M660828.\\
2020 $\mathnormal{Mathematics} \;\mathnormal{Subject}\;\mathnormal{Classification}$. 11F32, 11G55, 06A06.\\
$\mathnormal{Keywords:}$  Multiple zeta values, polylogarithms, partial orders, iterated integrals. }

\maketitle

\section{Introduction}\label{int}
    One variable multiple polylogarithms are iterated integrals of the differentials $\frac{dt}{t}$ and $\frac{dt}{1-t}$ from $0$ to $z$ (assume that $0<z<1$). For $k_1,k_2,\cdots,k_r\geq 1$, it is defined by
\[
\mathrm{Li}_{k_1,k_2,\cdots,k_r}(z)=\mathop{\int\cdots\int}_{0<t_1<t_2<\cdots<t_N<z}\omega_1(t_1)\omega_2(t_2)\cdots \omega_N(t_N), 0<z<1,
\]
where $N=k_1+k_2+\cdots+k_r$ and $\omega_i(t)=\frac{dt}{1-t}$ if $i\in\{1,k_1+1,\cdots,k_1+k_2+\cdots+k_{r-1}+1\}$, $\omega_i(t)=\frac{dt}{t}$ if $i\notin\{1,k_1+1,\cdots,k_1+k_2+\cdots+k_{r-1}+1\}$.

 For $k_1,\cdots, k_r\geq 1$, it has a series expression
\[
\mathrm{Li}_{k_1,k_2,\cdots,k_r}(z)=\sum_{0<n_1<n_2<\cdots<n_r}\frac{z^{n_r}}{n_1^{k_1}n_2^{k_2}\cdots n_r^{k_r}},\; |z|<1.
\]
When $r=1$ and $k_1=1$, $\mathrm{Li}_1(z)$ is the logarithm function
\[
\mathrm{Li}_1(z)=\sum_{n>0}\frac{z^n}{n}=-\mathrm{log}\; (1-z),\;|z|<1.
\]
When $k_r\geq 2$, the evaluation of $\mathrm{Li}_{k_1,k_2,\cdots,k_r}(z)$ at the point $z=1$ gives the classical multiple zeta value
\[
\mathrm{Li}_{k_1,k_2,\cdots,k_r}(z)|_{z=1}=\zeta({k_1,k_2,\cdots,k_r})=\sum_{0<n_1<n_2<\cdots<n_r}\frac{1}{n_1^{k_1}n_2^{k_2}\cdots n_r^{k_r}}.
\]

From Kontsevich's iterated integral representations of multiple zeta values, we know that  iterated integrals of the differentials $\frac{dt}{t}$ and $\frac{dt}{1-t}$ from $0$ to $1$ are multiple zeta values if they are convergent.
In this paper, we are interested in  iterated integrals on products of one variable multiple polylogarithms. 

Denote by ${\bf{P}}_{\mathbb{Q}}^{\,\mathrm{log}}$ the $\mathbb{Q}$-algebra which is generated by 
\[
1, \mathrm{Li}_{k_1,k_2,\cdots,k_r}(z),\;\mathrm{Li}_{k_1,k_2,\cdots,k_r}(1-z)
\]
for all $k_1,k_2,\cdots,k_r\geq 1$ and $r\geq 1$. Clearly all the functions in ${\bf{P}}_{\mathbb{Q}}^{\,\mathrm{log}}$ are convergent on  $0<z<1$. 
\begin{Thm}\label{main}
For $f_1(z),f_2(z),\cdots,f_N(z)\in {\bf{P}}_{\mathbb{Q}}^{\,\mathrm{log}}$, $\omega_1(z),\omega_2(z),\cdots, \omega_N(z)\in \{ \frac{dz}{z},\frac{dz}{1-z} \}$, if the iterated integral 
\[
I(f_1,\cdots,f_N;\omega_1,\cdots,\omega_N)=\mathop{\int\cdots\int}\limits_{0<z_1<z_2<\cdots<z_N<1}f_1(z_1)\omega_1(z_1)f_2(z_2)\omega_{2}(z_2)\cdots f_N(z_N)\omega_{N}(z_N)
\]
is convergent, then the value $I(f_1,\cdots,f_N;\omega_1,\cdots,\omega_N)$ is a $\mathbb{Q}$-linear combination of multiple zeta values.
\end{Thm}

\begin{rem}
About the convergence of the above iterated integrals,  from the proof of Theorem \ref{main}, it is easy to see that the integral
\[
I(f_1,\cdots,f_N;\omega_1,\cdots,\omega_N)
\] is convergent if and only 
if 
\[
\begin{split}
&f_1(z)\omega_1(z)\in {\bf{P}}_{\mathbb{Q}}^{\,\mathrm{log}}\cdot\frac{dz}{1-z}+\sum_{k_1,k_2,\cdots,k_r\geq 1}\mathrm{Li}_{k_1,k_2,\cdots,k_r}(z){\bf{P}}_{\mathbb{Q}}^{\,\mathrm{log}}\cdot\frac{dz}{z}\\
\end{split}
\]
and 
\[
f_N(z)\omega_N(z)\in {\bf{P}}_{\mathbb{Q}}^{\,\mathrm{log}}\cdot\frac{dz}{z}+\sum_{l_1,l_2,\cdots,l_s\geq 1}\mathrm{Li}_{l_1,l_2,\cdots,l_s}(1-z){\bf{P}}_{\mathbb{Q}}^{\,\mathrm{log}}\cdot\frac{dz}{1-z}
\]
 for $N\geq 2$ and  
\[
f_1(z)\omega_1(z)\in \sum_{{{k_1,k_2,\cdots,k_r\geq 1}}}\left(\mathrm{Li}_{k_1,k_2,\cdots,k_r}(z){\bf{P}}_{\mathbb{Q}}^{\,\mathrm{log}}\cdot\frac{dz}{z}+\mathrm{Li}_{k_1,k_2,\cdots,k_r}(1-z){\bf{P}}_{\mathbb{Q}}^{\,\mathrm{log}}\cdot\frac{dz}{1-z} \right)
\] for $N=1$.
\end{rem}

From the iterated integral representations of multiple zeta values, it is obvious that Theorem \ref{main} is true if $f_i=1$ for all $i$. It is quite interesting that $$I(f_1,\cdots,f_N;\omega_1,\cdots,\omega_N)$$ is still in the algebra of multiple zeta values in general cases.

If the iterated integral $I(f_1,\cdots,f_N;\omega_1,\cdots,\omega_N)$ is convergent, it also has a series expression. We will give  explicit calculation in some cases. But the strategy of  proving Theorem \ref{main} is not to compute the series expressions of the iterated integrals.

 Yamamoto \cite{yama} defined multiple integrals on finite partially ordered sets, and classical multiple zeta values can be viewed as multiple integrals on totally ordered sets. Furthermore, he proved that these multiple integrals on finite partially ordered sets are actually  $\mathbb{Q}$-linear combinations of classical multiple zeta values.

We will prove Theorem \ref{main} by  understanding the combinatorics of the iterated integrals. By combining this with the iterated integral expressions of one variable multiple polylogarithms, we will find that  $$I(f_1,\cdots,f_N;\omega_1,\cdots,\omega_N)$$ can be viewed as a multiple integral which is defined on a finite partially ordered set for any $f_1,\cdots,f_N\in {\bf{P}}_{\mathbb{Q}}^{\,\mathrm{log}}$. Since every partial order on  a finite set admits  refinements by some total orders, we get Theorem \ref{main}.

When the iterated integral $I(f_1,\cdots,f_N;\omega_1,\cdots,\omega_N)$ is divergent, we use the same technique in the theory of regularization of multiple zeta values
to define regularized iterated integrals from $0$ to $1$.
\begin{Thm}\label{reg}
For any $f_1(z),f_2(z),\cdots,f_N(z)\in {\bf{P}}_{\mathbb{Q}}^{\,\mathrm{log}}$ and any $\omega_1(z),\omega_2(z),\cdots, \omega_N(z)\in \{ \frac{dz}{z},\frac{dz}{1-z} \}$, we can define a regularized iterated integral
\[
I^{\,\mathrm{reg}}(f_1,\cdots,f_N;\omega_1,\cdots,\omega_N),
\]
which also satisfies the shuffle product of iterated integrals.
If the iterated integral $I(f_1,\cdots,f_N;\omega_1,\cdots,\omega_N)$ is convergent, we have \[
I^{\,\mathrm{reg}}(f_1,\cdots,f_N;\omega_1,\cdots,\omega_N)=I(f_1,\cdots,f_N;\omega_1,\cdots,\omega_N).
\]
Furthermore, $I^{\,\mathrm{reg}}(f_1,\cdots,f_N;\omega_1,\cdots,\omega_N)$ always belongs to the algebra of multiple zeta values.
\end{Thm}

In the last section, we calculate the  series expressions of some special iterated integrals on products of one variable multiple polylogarithms. Although they are still multiple zeta values, their series representations are quite different. As an application, we have
\begin{Thm}\label{app}
$(i)$(Seki and Yamamoto) $\zeta(3)=\sum\limits_{n_1,n_2\geq 1}\frac{1}{n_1n_2}\cdot\frac{1\;\;\cdot\;\; 2\;\;\cdot\;\;\cdots\;\;\cdot\;\; n_2}{n_1\cdot(n_1+1)\cdot\;\,\cdots\;\,\cdot  (n_1+n_2)}.$\\
 For $k,l\geq 1$,
 \[
\sum_{n_1,n_2\geq 1}\frac{1}{n_1^kn_2^l}\cdot\frac{1\;\;\cdot\;\; 2\;\;\cdot\;\;\cdots\;\;\cdot\;\; n_2}{n_1\cdot(n_1+1)\cdot\;\,\cdots\;\,\cdot  (n_1+n_2)}=
\begin{cases}
\zeta(k+1,\underbrace{1,1,\cdots,1}_{l-2},2)&\; l\geq2\\
\zeta(k+2)&\; l=1,
\end{cases}
\]
\[
\zeta(k+l+1)=\sum_{n_1\geq 1,0<m_1<m_2<\cdots<m_l}\frac{1}{n_1^km_1m_2\cdots m_l}\cdot \frac{1\;\;\cdot\;\; 2\;\;\cdot\;\;\cdots\;\;\cdot\;\; m_l}{n_1\cdot(n_1+1)\cdot\;\,\cdots\;\,\cdot  (n_1+m_l)}.
\]
$(ii)$ For $k_1,k_2,\cdots,k_r\geq 1$, the series 
\[
\sum_{n_1,n_2,\cdots,n_r\geq 1}\frac{1}{n_1^{k_1}n_2^{k_2}\cdots n_r^{k_r}}\cdot \frac{1}{n_1(n_1+n_2)\cdots (n_1+n_2+\cdots+n_r)},
\]
\[
\sum_{\substack{n_1,n_2,\cdots,n_r\geq 1\\n_i<m_i,1\leq i\leq r}}\frac{1}{n_1^{k_1}n_2^{k_2}\cdots n_r^{k_r}}\cdot \frac{1}{m_1(m_1+m_2)\cdots (m_1+m_2+\cdots+m_{r-1})(m_1+m_2+\cdots+m_r)^2}
\]
are convergent and their values are in the algebra of multiple zeta values.
\end{Thm}
We will prove Theorem \ref{app} by calculating some special iterated integrals on products of one variable multiple polylogarithms. In the proof, we also need some basic properties about the Beta function $B(\alpha,\beta)$ which is defined by
\[
B(\alpha,\beta)=\mathop{\int}_0^1t^{\alpha-1}(1-t)^{\beta-1}dt.
\]

Furthermore we show that every multiple zeta value has series representations like (i) and (ii) in the last section.
 More precisely, we have
\begin{Thm}\label{conv}(Seki and Yamamoto)
For every multiple zeta value $\zeta(k_1,k_2,\cdots,k_N)$ with $N>2$, there are $N-2$ ways to write $\zeta(k_1,k_2,\cdots,k_N)$ as 
\[
\begin{split}
&\;\;\;\;\;\;\zeta(k_1,k_2,\cdots,k_r)\\
&=\sum_{\substack{0<m_1<m_2<\cdots<m_s \\0<n_1<n_2<\cdots<n_{s^\prime}}}\frac{1}{m_1^{l_1}m_2^{l_2}\cdots m_s^{l_s}}\cdot\frac{1}{n_1^{j_1}n_2^{j_2}\cdots n_{s^{\prime}}^{j_{s^{\prime}}}}\cdot B(m_s+i,n_{s^{\prime}}+1-i)\\
&=\sum_{\substack{0<m_1<m_2<\cdots<m_s \\0<n_1<n_2<\cdots<n_{s^\prime}}}\frac{1}{m_1^{l_1}m_2^{l_2}\cdots m_s^{l_s}}\cdot \frac{1}{n_1^{j_1}n_2^{j_2}\cdots n_{s^{\prime}}^{j_{s^{\prime}}}}\\
&\;\;\;\;\;\;\;\;\;\;\;\;\;\;\;\;\;\;\;\;\;\;\;\;\;\;\;\;\;\;\;\;\;\;\;\;\;\cdot \frac{1\;\;\cdot\;\; 2\;\;\cdot\;\;\cdots\;\;\cdot\;\; (m_s+i-1)}{(n_{s^{\prime}}-i+1)\cdot(n_{s^{\prime}}-i+2)\cdot\;\,\cdots\;\,\cdot  (n_{s^{\prime}}+m_s)}\\
\end{split}
\]
for some $(l_1,l_2,\cdots,l_s), (j_1,j_2,\cdots,j_{s^{\prime}})$ and $ i$ satisfying
\[
l_1+l_2+\cdots+l_s+j_1+j_2+\cdots+j_{s^{\prime}}+1=N, i\in \{0,1\}.
\]

\end{Thm}

In fact, the above series representations for multiple zeta values are the connected sums defined by Seki and Yamamoto \cite{SY} to study the duality for multiple zeta values.
Thus Theorem \ref{app} (i) and Theorem \ref{conv} are already known to Seki and Yamamoto. But the approach here is different from \cite{SY}. All the connected sums can be viewed as iterated integrals on products of one variable multiple polylogarithms naturally.

For the calculation of iterated integrals on products of one variable multiple polylogarithms in general cases, we define  the multiple Beta functions and study their basic properties. We show that one can use the multiple Beta functions to calculate the iterated integrals in general cases by induction.

\section{Multiple integrals on finite partially ordered sets} 
In this section, we firstly give the definition of partial order on a (finite) set.
Secondly we  introduce the definition of multiple integrals on finite partially ordered sets, which is firstly defined by Yamamoto in \cite{yama}.  Lastly we  show that all multiple integrals on finite partially ordered sets are $\mathbb{Q}$-linear combinations of multiple zeta values.

\begin{Def}
A partial order $R$ on a (finite) set $V$ is a subset of $V^2=V\times V$ which satisfies:\\
(i) If $(t_1,t_2)\in R$ and $(t_2,t_3)\in R$, then $(t_1,t_3)\in R$;\\
(ii) For all $t\in V$, $(t,t)\notin R$.
\end{Def}

For a partial order $R$ on $V$, we usually write $t_1\prec t_2$ if $(t_1,t_2)\in R$.

Recall that multiple zeta value $\zeta(k_1,k_2,\cdots,k_r)$ where $k_r\geq 2$ can be defined by 
\[
\zeta(k_1,k_2,\cdots,k_r)=\mathop{\int\cdots\int}\limits_{0<t_1<t_2<\cdots<t_N<1}\omega_1(t_1)\omega_{2}(t_2)\cdots \omega_{N}(t_N),
\]
where $N=k_1+k_2+\cdots+k_r$, 
\[
\omega_i(t)=
\begin{cases}
\frac{dt}{1-t}, &\mathrm{if}\;i\in \{1,k_1+1,k_1+k_2+1,\cdots,k_1+\cdots+k_{r-1}+1 \},\\ 
\frac{dt}{t}, &\mathrm{else}.\\
\end{cases}
\]
It is clear that 
\[
t_1\prec t_2\prec \cdots \prec t_N
\]
actually defines a total order on the finite set $V=\{t_1,t_2,\cdots,t_N\}$. It corresponds to the complex 
\[
\{(t_1,t_2,\cdots,t_N)\mid 0<t_1<t_2<\cdots<t_N<1\} 
\]
in the cube $[0,1]^N$ naturally.

Inspired by this simple observation, Yamamoto \cite{yama} defined the multiple integrals on finite partially ordered sets.
\begin{Def}\label{gmzv}
For a finite set $V=\{t_1,t_2,\cdots, t_N\}$ and a partial order $R$ on $V$, define 
\[ 
\Delta_R=\{(t_1,t_2,\cdots,t_N)\mid \,0<t_i<t_j<1, \mathrm{iff}\; t_i\prec t_j \;\mathrm{in}\;R\}.
\]
For $\omega_i(t)\in \{\frac{dt}{t},\frac{dt}{1-t}\}$, if 
\[
\mathop{\int}\limits_{\Delta_R}\omega_1(t_1)\omega_{2}(t_2)\cdots \omega_{N}(t_N)
\]
is convergent, then it is denoted by $I^R(\omega_1,\cdots,\omega_N)$  and called the multiple integral on partially ordered  set $V$.
\end{Def}

The multiple integrals on totally ordered sets are surely multiple zeta values by definition.
Now we calculate some simple examples in general cases to give a rough impression about the above definition.

\begin{ex}\label{3}
For $V=\{t_1,t_2,t_3\},R=\{t_1\prec t_3,t_2\prec t_3\}, \omega_1(t)=\omega_2(t)=\frac{dt}{1-t}$ and $\omega_3(t)=\frac{dt}{t}$,
\[
\begin{split}
&\;\;\;\;\;I^R(\omega_1,\omega_2,\omega_3)\\
&=\mathop{\int}\limits_{{\substack {0<t_1<t_3<1\\0<t_2<t_3<1}}}\frac{dt_1}{1-t_1}\frac{dt_2}{1-t_2}\frac{dt_3}{t_3}  
=\int^1_0\left(\int_0^{t_3} \frac{dt}{1-t} \right)^2\frac{dt_3}{t_3}
=\sum_{n_1,n_2\geq 1}\frac{1}{n_1n_2(n_1+n_2)}.\\
\end{split}
\]
As a result,
\[
I^R(\omega_1,\omega_2,\omega_3)=\sum_{n_1,n_2\geq 1}\frac{1}{n_1^2}\left(\frac{1}{n_2}-\frac{1}{n_1+n_2} \right)=\sum_{n_1\geq 1}\frac{1}{n_1^2}\sum_{n_2=1}^{n_1}\frac{1}{n_2}=\zeta(3)+\zeta(1,2).
\]
Since $\zeta(1,2)=\zeta(3)$, we also have $I^R(\omega_1,\omega_2,\omega_3)=2\zeta(3)$.
\end{ex}

\begin{ex}\label{5}
For $V=\{t_1,t_2,t_3,t_4,t_5\}, R=\{t_1\prec t_2\prec t_3\prec t_5, t_1\prec t_4\prec t_5\}$, $\omega_1(t)=\omega_2(t)=\omega_3(t)=\frac{dt}{1-t}$ and $\omega_4(t)=\omega_5(t)=\frac{dt}{t}$.
\[
\begin{split}
&\;\;\;\;\;I^R(\omega_1,\omega_2,\omega_3,\omega_4,\omega_5)\\
&=\mathop{\int}\limits_{{\substack {0<t_1<t_2<t_3<t_5<1\\0<t_1<t_4<t_5<1}}}\frac{dt_1}{1-t_1}\frac{dt_2}{1-t_2}\frac{dt_3}{1-t_3}\frac{dt_4}{t_4}\frac{dt_5}{t_5}  
\\
\end{split}
\]
\[
\begin{split}
&=\mathop{\int}_{0<t_1<t_5<1}\left(\mathop{\int}_{t_1<t_2<t_3<t_5} \frac{dt_2}{1-t_2} \frac{dt_3}{1-t_3}\right)\left(\mathop{\int}_{t_1<t_4<t_5}\frac{dt_4}{t_4}\right)\frac{dt_1}{1-t_1}\frac{dt_5}{t_5}\\
&=\mathop{\int}_{0<t_1<t_5<1}\left(\mathop{\int}_{t_1<t_3<t_5}\left(\sum_{n_2\geq 1}\frac{t_3^{n_2}}{n_2}-\frac{t_1^{n_2}}{n_2}\right)\frac{dt_3}{1-t_3} \right)\left(\mathrm{log}\;t_5- \mathrm{log}\;t_1\right)\frac{dt_1}{1-t_1}\frac{dt_5}{t_5}\\
&=\mathop{\int}_{0<t_1<t_5<1}\sum_{n_2,n_3\geq 1}\left(\frac{t_5^{n_2+n_3}}{n_2(n_2+n_3)}-\frac{t_1^{n_2}t_5^{n_3}}{n_2n_3}-\frac{t_1^{n_2+n_3}}{n_2(n_2+n_3)}+\frac{t_1^{n_2+n_3}}{n_2n_3} \right)\mathrm{log}\;\frac{t_5}{t_1}\frac{dt_1}{1-t_1}\frac{dt_5}{t_5}.\\
\end{split}
\]
Since $\frac{d}{dt}\left( \frac{t^n\mathrm{log}\;t}{n}-\frac{t^n}{n^2}\right)=t^{n-1}\mathrm{log}\;t$ for $n\geq 1$, $$\int t^{n-1}\mathrm{log}\;tdt=\frac{t^{n}}{n}\left(\mathrm{log}\;t -\frac{1}{n}\right)+C,$$ we have 
\[
\begin{split}
&\;\;\;\;\;I^R(\omega_1,\omega_2,\omega_3,\omega_4,\omega_5)\\
&=\mathop{\int}_0^1\sum_{n_2,n_3,n_4\geq1}\left( \mathrm{log}\;t_5-\mathrm{log}\;t_1+\frac{1}{n_4} \right)\Bigg{[}\frac{t_5^{n_2+n_3}t_1^{n_4}}{n_2(n_2+n_3)n_4}-\frac{t_1^{n_2+n_4}t_5^{n_3}}{(n_2+n_4)n_2n_3}-\\
&\;\;\;\;\;\;\frac{t_1^{n_2+n_3+n_4}}{n_2(n_2+n_3)(n_2+n_3+n_4)}+\frac{t_1^{n_2+n_3+n_4}}{n_2n_3(n_2+n_3+n_4)}\Bigg{]}_{t_1=0}^{t_1=t_5}\frac{dt_5}{t_5}\\
\end{split}
\]
\[
\begin{split}
&=\mathop{\int}_0^1\sum_{n_2,n_3,n_4\geq 1}\Bigg{[}\frac{t_5^{n_2+n_3+n_4}}{n_2(n_2+n_3)n_4^2}-\frac{t_5^{n_2+n_3+n_4}}{(n_2+n_4)n_2n_3n_4}-\\
&\;\;\;\;\;\;\frac{t_5^{n_2+n_3+n_4}}{n_2(n_2+n_3)(n_2+n_3+n_4)n_4}+\frac{t_5^{n_2+n_3+n_4}}{n_2n_3(n_2+n_3+n_4)n_4}\Bigg{]}\frac{dt_5}{t_5}\\
&=\sum_{n_2,n_3,n_4\geq 1}\Bigg{[}\frac{1}{n_2(n_2+n_3)n_4^2(n_2+n_3+n_4)}-\frac{1}{(n_2+n_4)n_2n_3n_4(n_2+n_3+n_4)}-\\
&\;\;\;\;\;\;\frac{1}{n_2(n_2+n_3)(n_2+n_3+n_4)^2n_4}+\frac{1}{n_2n_3(n_2+n_3+n_4)^2n_4}\Bigg{]}.\\
\end{split}
\]
\end{ex}

The procedure to calculate the multiple integrals on finite partially ordered sets is very clear. Firstly we calculate the iterated integrals on some totally ordered subsets. Then we patch  all the data together and  reduce the multiple integrals to iterated integrals of power series in $t$ or $1-t$ which are convergent if $0<t<1$.

In Example \ref{3}, the multiple integral is still in the algebra of multiple zeta values. But in Example \ref{5}, it is not clear whether it is  a multiple zeta value or not at first sight. The interesting thing is that it is a $\mathbb{Q}$-linear combination of multiple zeta values. The following theorem has been already proved by Yamamoto in \cite{yama}. We give a different proof here . 
\begin{Thm}\label{gmv}
For any finite set $V$ and any partial order $R$ on $V$, the multiple integrals  on the partially ordered set $V$ are sums of finite numbers of multiple zeta values, if they are convergent.
\end{Thm}
 \noindent{\bf Proof}:
Let $V=\{t_1,t_2,\cdots,t_N\}, \omega_i(t)\in\{\frac{dt}{t},\frac{dt}{1-t}\}$ for $i=1,2,\cdots,N$, and assume that the integral $I^R(\omega_1,\omega_2,\cdots, \omega_N)$ is convergent. 

Recall that  $\Delta_R=\{(t_1,t_2,\cdots,t_N)\mid \,0<t_i<t_j<1, \mathrm{iff}\; t_i\prec t_j \;\mathrm{in}\;R\}$. 
Denote by $S_n$ the  group of permutations on the set $\{1,2,\cdots, n\}$.
Since the set of real numbers $\mathbb{R}$ is a totally ordered (infinite) set with the natural order $<$, the complex $\Delta_R$ can be written as 
\[
\Delta_R=\cup_{\sigma\in S_R}\{(t_1,t_2,\cdots,t_N)\mid 0<t_{\sigma(1)}<t_{\sigma(2)}<\cdots<t_{\sigma(N)}<1\}\bigcup L_R,
\]
where $$ S_R=\{\sigma\in S_n\mid \sigma(m)< \sigma(n) \;\mathrm{if}\; t_m\prec t_n\; \mathrm{in}\;R\}$$ and $\mathrm{dim}\;L_R<N$.
As a result of the above decomposition of the complex $\Delta_R$, we have 
\[
\begin{split}
&\;\;\;\;\;I^R(\omega_1,\omega_2,\cdots,\omega_N)\\
&= \mathop{\int}_{\Delta_R}\omega_1(t_1)\omega_2(t_2)\cdots \omega_N(t_N)=\sum_{\sigma\in S_R}\mathop{\int\cdots \int}_{0<t_{\sigma(1)}<t_{\sigma(2)}<\cdots<t_{\sigma(N)}<1}\omega_1(t_1)\omega_2(t_2)\cdots \omega_N(t_N)\\
&=\sum_{\sigma\in S_R}\mathop{\int\cdots \int}_{0<t_{1}<t_{2}<\cdots<t_{N}<1}\omega_1(t_{\sigma^{-1}(1)})\omega_2(t_{\sigma^{-1}(2)})\cdots \omega_N(t_{\sigma^{-1}(N)}).
\end{split}
\]
Since $\omega_i(t)=f_i(t)\,dt,1\leq i\leq N$ and $f_i(t)>0$ for $0<t<1$, we have that the integral $I^R(\omega_1,\omega_2,\cdots,\omega_N)$ is convergent if and only if the iterated integral
\[
\mathop{\int\cdots \int}_{0<t_{1}<t_{2}<\cdots<t_{N}<1}\omega_1(t_{\sigma^{-1}(1)})\omega_2(t_{\sigma^{-1}(2)})\cdots \omega_N(t_{\sigma^{-1}(N)})
\]
is convergent for all $ i=1,\cdots,k$.

Since $\omega_i(t)\in\{\frac{dt}{t},\frac{dt}{1-t}\}$, from the theory of multiple zeta values, we know that 
\[
\mathop{\int\cdots \int}_{0<t_{1}<t_{2}<\cdots<t_{N}<1}\omega_1(t_{\sigma^{-1}(1)})\omega_2(t_{\sigma^{-1}(2)})\cdots \omega_N(t_{\sigma^{-1}(N)})
\] 
is convergent if and only if $\omega_{\sigma(1)}(t)=\frac{dt}{1-t}$ and $ \omega_{\sigma(N)}(t)=\frac{dt}{t}$. Furthermore, when it is convergent, it is  a multiple zeta value by direct calculation of the iterated integral.

In conclusion, all  multiple integrals on finite partially ordered sets are sums of finite numbers of multiple zeta values.   $\hfill\Box$\\

From Theorem \ref{gmv}, the definition of multiple integrals on finite partially ordered sets  does not produce any new objects other than classical multiple zeta values. But from Example \ref{3} and Example \ref{5}, we know that it does give new integral representations and new series representations of multiple zeta values.

We hope that multiple integrals on finite partially ordered sets will be useful to investigate the relations among multiple zeta values and to understand the irrationality of multiple zeta values.

Now we are ready to prove Theorem \ref{main} in the introduction.\\
\noindent{\bf Proof of Theorem \ref{main}}: For one variable multiple polylogarithm
\[
\mathrm{Li}_{k_1,k_2,\cdots,k_r}(z)=\sum_{0<n_1<n_2<\cdots<n_r}\frac{z^{n_r}}{n_1^{k_1}n_2^{k_2}\cdots n_r^{k_r}},\; |z|<1,
\]
in this paper we only consider the positive real case ($z$ is a positive real number), it can be written as
\[
\mathrm{Li}_{k_1,k_2,\cdots,k_r}(z)=\mathop{\int\cdots\int}_{0<t_1<t_2<\cdots<t_n<z}\omega_1(t_1)\omega_2(t_2)\cdots \omega_n(t_n),\;0<z<1, \eqno{(1)}
\] 
where $n=k_1+k_2+\cdots+k_r$, $\omega_1(t)=\omega_{k_1+1}(t)=\cdots =\omega_{k_1+k_2+\cdots k_{r-1}+1}(t)=\frac{dt}{1-t}$ and $\omega_i(t)=\frac{dt}{t}$ for other $i$.

Replacing $z$ by $1-z$ in $(1)$, we have 
\[
\begin{split}
&\;\;\;\;\;\;\;\mathrm{Li}_{k_1,k_2,\cdots,k_r}(1-z)\\
&=\mathop{\int\cdots\int}_{0<t_1<t_2<\cdots<t_n<1-z}\omega_1(t_1)\omega_2(t_2)\cdots \omega_n(t_n),\\
&=(-1)^n\mathop{\int\cdots\int}_{1>t_1>t_2>\cdots>t_n>z}\omega_1(1-t_1)\omega_2(1-t_2)\cdots \omega_n(1-t_n)\\
&=\mathop{\int\cdots\int}_{z<t_n<\cdots<t_2<t_1<1}\omega^{\prime}_n(t_n)\omega^{\prime}_{n-1}(t_{n-1})\cdots \omega^{\prime}_1(t_1),\\
\end{split}
\]
where $\,\omega_i^{\prime}(t)+\omega_i(t)=\frac{dt}{t(1-t)}$, for $ i=1,2,\cdots,n$. Thus
\[
\mathrm{Li}_{k_1,k_2,\cdots,k_r}(1-z)=\mathop{\int\cdots\int}_{z<t_n<\cdots<t_2<t_1<1}\omega^{\prime}_n(t_n)\omega^{\prime}_{n-1}(t_{n-1})\cdots \omega^{\prime}_1(t_1). \eqno{(2)}
\]

For $f_1(z),f_2(z),\cdots,f_N(z)\in {\bf{P}}_{\mathbb{Q}}^{\,\mathrm{log}}$ and $\omega_1(z),\omega_2(z),\cdots, \omega_N(z)\in \{ \frac{dz}{z},\frac{dz}{1-z} \}$, in order to prove that the (convergent) iterated  integral 
\[
I(f_1,\cdots,f_N;\omega_1,\cdots,\omega_N)=\mathop{\int\cdots\int}\limits_{0<z_1<z_2<\cdots<z_N<1}f_1(z_1)\omega_1(z_1)f_2(z_2)\omega_{2}(z_2)\cdots f_N(z_N)\omega_{N}(z_N)
\]
is actually a $\mathbb{Q}$-linear combination of multiple zeta values, without loss of generality we can assume that $f_1(z),f_2(z),\cdots,f_N(z)$ are all monic monomials in 
\[
1,\mathrm{Li}_{k_1,k_2,\cdots,k_r}(z),\mathrm{Li}_{k_1,k_2,\cdots,k_r}(1-z), k_1,\cdots,k_r\geq 1,r\geq1.
\]

Assume that $$f_i(z)=\prod\limits_{k_1,\cdots,k_r\geq 1} \left(\mathrm{Li}_{k_1,k_2,\cdots,k_r}(z)\right)^{m(f_i)_{k_1,k_2,\cdots,k_r}}\cdot \prod\limits_{k_1,\cdots,k_r\geq 1} \left(\mathrm{Li}_{k_1,k_2,\cdots,k_r}(1-z)\right)^{m^{\prime}(f_i)_{k_1,k_2,\cdots,k_r}},$$ where $m(f_i)_{k_1,k_2,\cdots,k_r},m^{\prime}(f_i)_{k_1,k_2,\cdots,k_r}>0$, for $1\leq i\leq N$.

Let 
\[
K_i=\sum_{k_1,\cdots,k_r\geq 1}\Big{[} m(f_i)_{k_1,k_2,\cdots,k_r}(k_1+k_2+\cdots+k_r)+ m^{\prime}(f_i)_{k_1,k_2,\cdots,k_r}(k_1+k_2+\cdots+k_r)\Big{]},
\]
$$K=K_1+K_2+\cdots+K_N,T=K+N,$$
define a partial order $R$ on the set \[
\begin{split}
&V_T=\\
&\Bigg{\{}{}_it_1^{k_1,k_2,\cdots,k_r}, {}_it_2^{k_1,k_2,\cdots,k_r},\cdots,{}_it_{m(f_i)_{k_1,k_2,\cdots,k_r}(k_1+k_2+\cdots+k_r)}^{k_1,k_2,\cdots,k_r},{}_is_1^{k_1,k_2,\cdots,k_r}, {}_is_2^{k_1,k_2,\cdots,k_r},  \\
& \cdots,{}_is_{m^{\prime}(f_i)_{k_1,k_2,\cdots,k_r}(k_1+k_2+\cdots+k_r)}^{k_1,k_2,\cdots+\,k_r},\;z_i\,\Bigg{|}\,1\leq i\leq N, k_1,k_2,\cdots,k_r\geq 1\Bigg{\}}         \\
\end{split}
\]
 by
\[
\begin{split}
&{}_it_1^{k_1,k_2,\cdots,k_r}\prec {}_it_2^{k_1,k_2,\cdots,k_r}\prec \cdots \prec {}_it_{k_1+k_2+\cdots k_r}^{k_1,k_2,\cdots,k_r}\prec z_i,\\
&{}_1t_{k_1+k_2+\cdots k_r+1}^{k_1,k_2,\cdots,k_r}\prec {}_it^{k_1,k_2,\cdots,k_r}_{k_1+k_2+\cdots +k_r+2}\prec\cdots \prec {}_it_{2(k_1+k_2+\cdots+k_r)}^{k_1,k_2,\cdots,k_r}\prec z_i,\\
&\cdots\;\;\;\;\;\;\;\;\;\;\;\;\cdots\;\;\;\;\;\;\;\;\;\;\;\;\cdots\\
&{}_it^{k_1,k_2,\cdots,k_r}_{(m(f_i)_{k_1,k_2,\cdots,k_r}-1)(k_1+k_2+\cdots+k_r)+1}\prec{}_it^{k_1,k_2,\cdots,k_r}_{(m(f_i)_{k_1,k_2,\cdots,k_r}-1)(k_1+k_2+\cdots+k_r)+2}\prec\cdots \\
&\;\;\;\;\;\;\;\;\;\;\;\;\;\;\;\;\;\;\;\;\;\;\;\;\;\;\;\;\prec {}_it^{k_1,k_2,\cdots,k_r}_{m(f_i)_{k_1,k_2,\cdots,k_r}(k_1+k_2+\cdots+k_r)}\prec z_i,
\end{split}
\]
\[
\begin{split}
z_i\prec {}_is_{k_1+k_2+\cdots+k_r}^{k_1,k_2,\cdots,k_r}\prec {}_is^{k_1,k_2,\cdots,k_r}_{k_1+k_2+\cdots+k_r-1}\prec \cdots \prec {}_is_1^{k_1,k_2,\cdots,k_r},
\end{split}
\]
\[
\begin{split}
&z_i\prec {}_is_{2(k_1+k_2+\cdots+k_r)}^{k_1,k_2,\cdots,k_r}\prec {}_is_{2(k_1+k_2+\cdots+k_r)-1}^{k_1,k_2,\cdots,k_r}\prec \cdots \prec {}_is_{k_1+k_2+\cdots+k_r+1}^{k_1,k_2,\cdots,k_r},\\
&\cdots\;\;\;\;\;\;\;\;\;\;\;\;\cdots \;\;\;\;\;\;\;\;\;\;\;\;\cdots \\
&z_i\prec {}_is_{m(f_i)_{k_1,k_2,\cdots,k_r}(k_1+k_2+\cdots+k_r)}^{k_1,k_2,\cdots,k_r}\prec {}_is_{m(f_i)_{k_1,k_2,\cdots,k_r}(k_1+k_2+\cdots+k_r)-1}^{k_1,k_2,\cdots,k_r}\prec\cdots\\
& \;\;\;\;\;\;\;\;\;\;\;\;\;\;\;\;\;\;\;\;\;\;\;\;\;\;\;\;\prec{}_is_{(m(f_i)_{k_1,k_2,\cdots,k_r}-1)(k_1+k_2+\cdots+k_r)+1}^{k_1,k_2,\cdots,k_r}
\end{split}
\]
\[
z_1\prec z_2 \prec \cdots \prec z_N, 
\]
for $i=1,2,\cdots,N$ and $k_1,k_2,\cdots,k_r\geq 1$.

From $(1)$ and $(2)$, it is clear that
\[
I(f_1,\cdots,f_N; \omega_1,\cdots,\omega_N)=\mathop{\int}_{\Delta_R}\omega_1\cdots \omega_T,
\]
where $\omega_1,\omega_2,\cdots,\omega_T$ are  differentials in  variables
\[
{}_it_1^{k_1,k_2,\cdots,k_r}, {}_it_2^{k_1,k_2,\cdots,k_r},\cdots,{}_it_{m(f_i)_{k_1,k_2,\cdots,k_r}(k_1+k_2+\cdots+k_r)}^{k_1,k_2,\cdots,k_r},
\]
\[
{}_is_1^{k_1,k_2,\cdots,k_r}, {}_is_2^{k_1,k_2,\cdots,k_r}, 
\cdots,{}_is_{m^{\prime}(f_i)_{k_1,k_2,\cdots,k_r}(k_1+k_2+\cdots+k_r)}^{k_1,k_2,\cdots+\,k_r}, z_i, 
\]
for $1\leq i\leq N$ and $k_1,k_2,\cdots,k_r\geq 1$ and $\omega_i(t)\in \{\frac{dt}{t},\frac{dt}{1-t}\}$ for $i\in \{1,\cdots,N\}$. Thus Theorem \ref{main} follows from Theorem \ref{gmv}. $\hfill\Box$\\

From the above proof, we know that 
\begin{Cor}\label{rea}
Every iterated integral on products of one variable multiple polylogarithms can be realized as a multiple integral on a finite partially ordered set in a natural way.
\end{Cor}

\section{Regularized iterated integral}

In this section, we firstly review the theory of shuffle regularization of multiple zeta values (see \cite{deligne}, \cite{ihara}, \cite{LM}).  By essentially the same method, we define regularized multiple integrals on finite partially ordered sets. Iterated integrals on products of one variable multiple polylogarithms can be viewed as multiple integrals on finite partially ordered sets in a natural way. Thus  we define the regularized iterated integral.

Denote by $\mathfrak{h}=(\mathbb{Q}\langle e^0,e^1\rangle,\rotatebox{90}{$\rotatebox{180}{$\exists$}$}\,)$, where $\mathbb{Q}\langle e^0,e^1\rangle$ means the  non-commutative plolynomial ring in variables $e^0,e^1$ and $\rotatebox{90}{$\rotatebox{180}{$\exists$}$}$ is the shuffle product defined by 
\[
w\;\rotatebox{90}{$\rotatebox{180}{$\exists$}$}\; 1=1\;\rotatebox{90}{$\rotatebox{180}{$\exists$}$}\;w=w,
\]
\[
xu\;\rotatebox{90}{$\rotatebox{180}{$\exists$}$}\;yv=x(u\;\rotatebox{90}{$\rotatebox{180}{$\exists$}$}\;yv )+y(u\;\rotatebox{90}{$\rotatebox{180}{$\exists$}$}\;yv ), x,y\in \{e^0,e^1\},
\]
inductively.

For $\epsilon>0,\eta>0$ and $\epsilon+\eta<1$, define a map $I^{(\epsilon,\eta)}:\mathfrak{h}\rightarrow \mathbb{R}$ by
\[
I^{(\epsilon,\eta)}(u_1u_2\cdots u_k)=\mathop{\int\cdots\int}_{\epsilon<t_1<t_2<\cdots<t_k<1-\eta}\omega_{u_1}(t_1)\omega_{u_2}(t_2)\cdots \omega_{u_k}(t_k), 
\]
for $u_1,u_2,\cdots,u_k\in \{e^0,e^1\}$ and $\omega_{e^0}(t)=\frac{dt}{t},\omega_{e^1}(t)=\frac{dt}{1-t}$. It is easy to check that the limit 
\[
\mathop{\mathrm{lim}}\limits_{\epsilon,\eta\rightarrow 0^+} I^{(\epsilon,\eta)}(u_1u_2\cdots u_k)
\]
exists if and only if $u_1=e^1$ and $u_k=e^0$. Moreover, the above limit is in the algebra of multiple zeta values if $u_1=e^1$ and $u_k=e^0$.

In general cases, we have 
\[
I^{(\epsilon,\eta)}(u_1u_2\cdots u_k)=P_{u_1,u_2,\cdots,u_k}(\mathrm{log}\;\epsilon,\mathrm{log}\;\eta)+O(\epsilon |\mathrm{log}\;\epsilon|^A +\eta|\mathrm{log}\;\eta|^B) \eqno{(3)}
\]
for some constants $A$ and $B$ (which depend on the word $u_1u_2\cdots u_k$), where $P_{u_1,u_2,\cdots,u_k}$ is a  polynomial of two variables with all coefficients in the algebra of multiple zeta values.
See the Appendix in \cite{LM} for a nice proof of formula $(3)$. We can also prove formula $(3)$ by direct calculation and induction.

\begin{Def}\label{mzv}
Define $I^{reg}:\mathfrak{h}\rightarrow \mathbb{R}$ by 
$I^{reg}(u_1u_2\cdots u_k)= P_{u_1,u_2,\cdots,u_k}(0,0)$.
\end{Def}

From the shuffle product on iterated integrals \cite{chen}, we know that $I^{reg}$ is a $\mathbb{Q}$-algebra homomorphism. Moreover,  images of elements in $\mathfrak{h}$ belong to the algebra of multiple zeta values.

For a finite set $V=\big{\{}t_1,t_2,\cdots,t_N\big{\}}$ and a partial order $R$ on $V$,  we define 
\[
I^{R,(\epsilon,\eta)}(\omega_1,\omega_2,\cdots,\omega_N)=\mathop{\int}_{\Delta_R\cap \;(\epsilon,1-\eta)^N}\omega_1\omega_2\cdots \omega_N.
\]
Then, by  the essentially same reason in the proof of Theorem \ref{gmv}, we have
\[
I^{R,(\epsilon,\eta)}(\omega_1,\omega_2,\cdots,\omega_N)=I^{(\epsilon,\eta)}(w)
\] 
for some $w\in \mathfrak{h}$.
Thus we have 
\[
I^{R,(\epsilon,\eta)}(\omega_1,\omega_2,\cdots,\omega_N)=P^R_{\omega_1,\omega_2,\cdots,\omega_k}(\mathrm{log}\;\epsilon,\mathrm{log}\;\eta)+O(\epsilon |\mathrm{log}\;\epsilon|^A +\eta|\mathrm{log}\;\eta|^B),
\]
where $P^R_{u_1,u_2,\cdots,u_k}$ is a polynomial  of two variables   with  coefficients in the algebra of multiple zeta values.
\begin{Def}\label{rmzv}
For $\omega_i\in \big{\{}\frac{dt}{t},\frac{dt}{1-t}\big{\}}$, define $I^{R,reg}(\omega_1,\omega_2,\cdots,\omega_N)=P^R_{\omega_1,\omega_2,\cdots,\omega_k}(0,0)$.
\end{Def}

Now we are ready to prove Theorem \ref{reg}.\\
\noindent{\bf Proof of Theorem \ref{reg}}: From the proof of Theorem \ref{main}, every iterated integral on products of one variable multiple polylogarithms (which may  be divergent) 
\[
I(f_1,f_2,\cdots,f_N;\omega_1,\omega_2,\cdots,\omega_N)
\]
(assume that $f_1,f_2,\cdots,f_N$ are all monic polynomials in one variable multiple polylogarithms) can be viewed as a multiple integral on a finite partially ordered set in a natural way. So we can use Definition \ref{rmzv} to define $$I^{reg}(f_1,f_2,\cdots,f_N;\omega_1,\omega_2,\cdots,\omega_N).$$

The shuffle product property of $I^{reg}(f_1,f_2,\cdots,f_N;\omega_1,\omega_2,\cdots,\omega_N)$ follows from the definition. From Theorem \ref{gmv} and formula $(3)$, we know that $$I^{reg}(f_1,f_2,\cdots,f_N;\omega_1,\omega_2,\cdots,\omega_N)$$ is still a $\mathbb{Q}$-linear combination of multiple zeta values.   $\hfill\Box$\\

\begin{rem}
In this paper, we only discuss the shuffle product on multiple zeta values. For another product which comes from the series representations of multiple zeta values, see \cite{hoff},\cite{ihara}.
\end{rem}

\section{Explicit calculation on iterated integrals}
Although iterated integrals on products of one variable multiple polylogarithms do not provide us anything which are different from classical multiple zeta values, it is still interesting to calculate their series expressions. We will calculate some interesting examples in this section.

For $r=1,f_1=\mathrm{Li}_1(z)\mathrm{Li}_1(1-z)$ and $\omega_1(z)=\frac{dz}{z}$,
\[
I(f_1;\omega_1)=\mathop{\int}_0^1\mathrm{Li}_1(z)\,\mathrm{Li}_1(1-z)\frac{dz}{z}=\mathop{\int}_0^1\sum_{n_1,n_2\geq 1}\frac{z^{n_1}(1-z)^{n_2}}{n_1n_2}\cdot\frac{dz}{z}=\sum_{n_1,n_2\geq 1}\frac{B(n_1,n_2+1)}{n_1n_2},
\]
where
\[
B(\alpha,\beta)=\int^1_0z^{\alpha-1}(1-z)^{\beta-1}dz=\frac{\Gamma{(\alpha)}\Gamma{(\beta)}}{\Gamma(\alpha+\beta)}
\]
is the Beta function (see \cite{rudin} for example).
Since $$B(n_1,n_2+1)=\frac{\Gamma(n_1)\Gamma(n_2+1)}{\Gamma(n_1+n_2+1)}=\frac{(n_1-1)!n_2!}{(n_1+n_2)!}=\frac{1\;\;\cdot\;\; 2\;\;\cdot\;\,\cdots\;\,\cdot\;\; n_2}{n_1\cdot(n_1+1)\cdot\;\,\cdots\;\,\cdot  (n_1+n_2)},$$ we have
\[
I(f_1;\omega_1)=\sum_{n_1,n_2\geq 1}\frac{1}{n_1n_2}\cdot \frac{(n_1-1)!n_2!}{(n_1+n_2)!}=\sum_{n_1,n_2\geq 1}\frac{1}{n_1n_2}\cdot\frac{1\;\;\cdot\;\; 2\;\;\cdot\;\;\cdots\;\;\cdot\;\; n_2}{n_1\cdot(n_1+1)\cdot\;\,\cdots\;\,\cdot  (n_1+n_2)}.
\]

On the other hand, let $V=\big{\{}t_1,t_2,z_1\big{\}}$, $R=\big{\{}t_1\prec z_1,z_1\prec t_2 \big{\}}$, $$\omega_1(t_1)=\frac{dt_1}{1-t_1},\omega_2(t_2)=\frac{dt_2}{t_2},\omega_3({z_1})=\frac{dz_1}{z_1},$$ thus
\[
I(f_1;\omega_1)=\mathop{\int}_{\Delta_R}\omega_1(t_1)\omega_2(t_2)\omega_3(z_1)=\mathop{\int\int\int}_{0<t_1<z_1<t_2<1}\frac{dt_1}{1-t_1}\frac{dz_1}{z_1}\frac{dt_2}{t_2}=\zeta(3).
\]

As an application, we prove that 
\[
\zeta(3)=\sum_{n_1,n_2\geq 1}\frac{1}{n_1n_2}\cdot\frac{1\;\;\cdot\;\; 2\;\;\cdot\;\;\cdots\;\;\cdot\;\; n_2}{n_1\cdot(n_1+1)\cdot\;\,\cdots\;\,\cdot  (n_1+n_2)}.\eqno{(4)}
\]
Formula $(4)$ is already known to many experts. It is quite interesting that formula $(4)$ can be viewed as two different realizations of one simple multiple integral on a finite partially ordered set.

Similarly, for $k,l\geq 2$,
\[
\begin{split}
&\;\;\;\;\mathop{\int}_0^1\mathrm{Li}_k(z)\mathrm{Li}_l(1-z)\frac{dz}{z}\\
&=\mathop{\int}_0^1\sum_{n_1,n_2\geq 1}\frac{z^{n_1}(1-z)^{n_2}}{n_1^kn_2^l}\cdot\frac{dz}{z}=\sum_{n_1,n_2\geq 1}\frac{1}{n_1^kn_2^l}\cdot\frac{1\;\;\cdot\;\; 2\;\;\cdot\;\;\cdots\;\;\cdot\;\; n_2}{n_1\cdot(n_1+1)\cdot\;\,\cdots\;\,\cdot  (n_1+n_2)},
\end{split}
\]
\[
\begin{split}
&\;\;\;\;\mathop{\int}_0^1\mathrm{Li}_k(z)\mathrm{Li}_l(1-z)\frac{dz}{z}\\
&=\mathop{\int\cdots \int}_{0<t_1<t_2<\cdots<t_k<z<s_l<\cdots<s_2<s_1<1}\frac{dt_1}{1-t_1}\frac{dt_2}{t_2}\cdots \frac{dt_k}{t_k}\frac{dz}{z}\frac{ds_k}{1-s_l}\cdots\frac{ds_2}{1-s_2}\frac{ds_1}{s_1}\\
&=\zeta(k+1,\underbrace{1,1,\cdots,1}_{l-2},2).   \\
\end{split}
\]
Thus for $k,l\geq 2$,
\[
\zeta(k+1,\underbrace{1,1,\cdots,1}_{l-2},2)=\sum_{n_1,n_2\geq 1}\frac{1}{n_1^kn_2^l}\cdot\frac{1\;\;\cdot\;\; 2\;\;\cdot\;\;\cdots\;\;\cdot\;\; n_2}{n_1\cdot(n_1+1)\cdot\;\,\cdots\;\,\cdot  (n_1+n_2)}. \eqno{(5)}
\]

By considering $\int_0^1\mathrm{Li}_k(z)\mathrm{Li}_1(1-z)\frac{dz}{z}$ for $k\geq 1$, we have 
\[
\zeta(k+2)=\sum_{n_1,n_2\geq 1}\frac{1}{n_1^k n_2}\cdot \frac{1\;\;\cdot\;\; 2\;\;\cdot\;\;\cdots\;\;\cdot\;\; n_2}{n_1\cdot(n_1+1)\cdot\;\,\cdots\;\,\cdot  (n_1+n_2)}. \eqno{(6)}
\]

For $r=1,f_1(z)=\mathrm{Li}_k(z)\;\mathrm{Li}_{\underbrace{\scriptstyle 1,1,\cdots,1}_{l}}(1-z)$ and  $\omega(z)=\frac{dz}{z}$, we have 
\[
\begin{split}
&\;\;\;\;\mathop{\int}_0^1\mathrm{Li}_k(z)\;\mathrm{Li}_{\underbrace{\scriptstyle 1,1,\cdots,1}_{l}}(1-z)\frac{dz}{z}\\
&=\mathop{\int}_0^1 \sum_{n_1\geq 1,0<m_1<m_2<\cdots<m_l}\frac{1}{n_1^km_1m_2\cdots m_l} z^{n_1}(1-z)^{m_l}\frac{dz}{z}\\
&=\sum_{n_1\geq 1,0<m_1<m_2<\cdots<m_l}\frac{1}{n_1^km_1m_2\cdots m_l} B(n_1,m_l)\\
&=\sum_{n_1\geq 1,0<m_1<m_2<\cdots<m_l}\frac{1}{n_1^km_1m_2\cdots m_l}\cdot \frac{1\;\;\cdot\;\; 2\;\;\cdot\;\;\cdots\;\;\cdot\;\; m_l}{n_1\cdot(n_1+1)\cdot\;\,\cdots\;\,\cdot  (n_1+m_l)}
\end{split}
\]
and 
\[
\begin{split}
&\;\;\;\;\mathop{\int}_0^1\mathrm{Li}_k(z)\;\mathrm{Li}_{\underbrace{\scriptstyle 1,1,\cdots,1}_{l}}(1-z)\frac{dz}{z}\\
&=\mathop{\int\cdots\int}_{0<t_1<t_2<\cdots<t_k<z<s_l<\cdots<s_{2}<s_1<1}\frac{dt_1}{1-t_1}\frac{dt_2}{t_2}\cdots \frac{dt_k}{t_k}\frac{dz}{z}\frac{ds_k}{s_k}\cdots \frac{ds_2}{s_2}\frac{ds_1}{s_1}  \\
&= \zeta{(k+l+1)}  . \\
\end{split}
\]
Thus 
\[
\zeta(k+l+1)=\sum_{n_1\geq 1,0<m_1<m_2<\cdots<m_l}\frac{1}{n_1^km_1m_2\cdots m_l}\cdot \frac{1\;\;\cdot\;\; 2\;\;\cdot\;\;\cdots\;\;\cdot\;\; m_l}{n_1\cdot(n_1+1)\cdot\;\,\cdots\;\,\cdot  (n_1+m_l)}. \eqno{(7)}
\]

\noindent{\bf Proof of Theorem \ref{conv}:}
In fact, for the multiple zeta value $\zeta(k_1,k_2,\cdots,k_r), N=k_1+k_2+\cdots+k_r$, we have
\[
\begin{split}
&\;\;\;\;\;\;\;\zeta(k_1,k_2,\cdots,k_r)\\
&=\mathop{\int\cdots \int}_{0<z_1<z_2<\cdots<z_N<1}\omega_1(z_1)\omega_2(z_2)\cdots \omega_N(t_N)\\
&=\mathop{\int}_0^1 \left(\mathop{\int\cdots\int}_{0<z_1<\cdots<z_{i-1}<z} \omega_1(z_1)\cdots \omega_{i-1}(z_{i-1})  \right)\left(\mathop{\int\cdots\int}_{z<z_{i+1}<\cdots<z_N<1}\omega_{i+1}(z_{i+1})\cdots \omega(z_N) \right)\omega_i(z)\\
&=\mathop{\int}_0^1\left(\mathop{\int\cdots\int}_{0<z_1<\cdots<z_{i-1}<z} \omega_1(z_1)\cdots \omega_{i-1}(z_{i-1})  \right)\left(\mathop{\int\cdots\int}_{1-z>z_{i+1}>\cdots>z_N>0}\omega^{\prime}_{i+1}(z_{i+1})\cdots \omega^{\prime}(z_N) \right)\omega_i(z)\\
&=\mathop{\int}_0^1 \mathrm{Li}_{l_1,l_2,\cdots,l_s}(z)\mathrm{Li}_{j_1,j_2,\cdots,j_{s^{\prime}}}(1-z)\omega_i(z)
\end{split}
\]
for some $\mathrm{Li}_{l_1,l_2,\cdots,l_s}(z),\mathrm{Li}_{j_1,j_2,\cdots,j_{s^{\prime}}}(1-z)$ and $l_s+j_{s^{\prime}}+1=N,1<i<N$. Thus every multiple zeta value has series representations like $(5)$ and $(6)$.   $\hfill\Box$\\

\begin{rem}
To study the duality of multiple zeta values, Seki and Yamamoto \cite{SY} defined the connected sum
\[
Z({\overrightarrow{k}},\overrightarrow{l})=\sum_{\substack{0<m_1<m_2<\cdots<m_r  \\0<n_1<n_2<\cdots<n_s  }}\frac{1}{m_1^{k_1}m_2^{k_2}\cdots m_r^{k_r}}\frac{1}{n_1^{l_1}n_2^{l_2}\cdots n_s^{l_s}}\cdot \frac{m_r!n_s!}{(m_r+n_s)!}
\]
for every $\overrightarrow{k}=(k_1,k_2,\cdots,k_r)$ and $\overrightarrow{l}=(l_1,l_2,\cdots,l_s)$.
Since 
\[
\begin{split}
&\;\;\;\;\;Z({\overrightarrow{k}},\overrightarrow{l})\\
&=\sum_{\substack{0<m_1<m_2<\cdots<m_r  \\0<n_1<n_2<\cdots<n_s  }}\frac{1}{m_1^{k_1}m_2^{k_2}\cdots m_r^{k_r-1}}\frac{1}{n_1^{l_1}n_2^{l_2}\cdots n_s^{l_s}}\cdot \frac{(m_r-1)!n_s!}{(m_r+n_s)!}   \\
&= \sum_{\substack{0<m_1<m_2<\cdots<m_r  \\0<n_1<n_2<\cdots<n_s  }}\frac{1}{m_1^{k_1}m_2^{k_2}\cdots m_r^{k_r}}\frac{1}{n_1^{l_1}n_2^{l_2}\cdots n_s^{l_s-1}}\cdot \frac{m_r!(n_s-1)!}{(m_r+n_s)!},  \\
\end{split}
\]
 the new series representations of every multiple zeta value in Theorem \ref{conv} are actually connected sums defined by Seki and Yamamoto in \cite{SY}.  But the proof of Theorem \ref{conv} here is very different from \cite{SY}.
\end{rem}

For $f_1(z)=\mathrm{Li}_{k_1}(z),f_2(z)=\mathrm{Li}_{k_2}(z),\cdots,f_r(z)=\mathrm{Li}_{k_r}(z)$ and $\omega_1(z)=\cdots=\omega_r(z)=\frac{dz}{z}$, we have
\[
\begin{split}
&\;\;\;\;I(f_1,\cdots,f_r;\omega_1,\cdots,\omega_r)\\
&=\mathop{\int\cdots\int}_{0<z_1<z_2<\cdots<z_r<1}\mathrm{Li}_{k_1}(z_1)\frac{dz_1}{z_1}\mathrm{Li}_{k_2}(z_2)\frac{dz_2}{z_2}\cdots \mathrm{Li}_{k_r}(z_r)\frac{dz_r}{z_r}\\
&=  \mathop{\int\cdots\int}_{0<z_2<\cdots<z_r<1}\sum_{n_1\geq 1}\frac{z_2^{n_1}}{n_1^{k_1+1}}\mathrm{Li}_{k_2}(z_2)\frac{dz_2}{z_2}\cdots \mathrm{Li}_{k_r}(z_r)\frac{dz_r}{z_r}\\
&=\;\cdots\\
&=\sum_{n_1,n_2,\cdots,n_r\geq 1}\frac{1}{n_1^{k_1}n_2^{k_2}\cdots n_r^{k_r}}\cdot \frac{1}{n_1(n_1+n_2)\cdots (n_1+n_2+\cdots+n_r)}.
\end{split}
\]
As a result of Theoren \ref{gmv}, we have
\begin{Cor}
For any $k_1,k_2,\cdots,k_r\geq 1$, the series 
\[
\sum_{n_1,n_2,\cdots,n_r\geq 1}\frac{1}{n_1^{k_1}n_2^{k_2}\cdots n_r^{k_r}}\cdot \frac{1}{n_1(n_1+n_2)\cdots (n_1+n_2+\cdots+n_r)}
\]
is convergent and its value is in the algebra of multiple zeta values.
\end{Cor}

Let $f_1(z)=\mathrm{Li}_{k_1}(z),f_2(z)=\mathrm{Li}_{k_2}(z),\cdots,f_r(z)=\mathrm{Li}_{k_r}(z)$,    $f_{r+1}(z)=1$, and $\omega_1(z)=\cdots=\omega_r(z)=\frac{dz}{1-z}$, $\omega_{r+1}(z)=\frac{dz}{z}$,  similarly we have
\begin{Cor}
For any $k_1,k_2,\cdots,k_r\geq 1$, the seires 
\[
\sum_{\substack{n_1,n_2,\cdots,n_r\geq 1\\n_i<m_i,1\leq i\leq r}}\frac{1}{n_1^{k_1}n_2^{k_2}\cdots n_r^{k_r}}\cdot \frac{1}{m_1(m_1+m_2)\cdots (m_1+m_2+\cdots+m_{r-1})(m_1+m_2+\cdots+m_r)^2}
\]
is convergent and its value is in the algebra of multiple zeta values.\\
\end{Cor}

Let $f_1(z)=\mathrm{Li}_{k_1}(z)\mathrm{Li}_{k_2}(z)\cdots \mathrm{Li}_{k_s}(z), f_2(z)=\cdots=f_r(z)=1$, and $ \omega_1(z)=\cdots=\omega_r(z)=\frac{dz}{z}$, we have 
\begin{Cor}\label{zho}
For $r,s \geq 1$ and $k_1,k_2,\cdots,k_s\geq 1$, the series 
\[
\sum_{n_1,n_2,\cdots,n_s\geq 1}\frac{1}{n_1^{k_1}n_2^{k_2}\cdots n_s^{k_s}(n_1+n_2+\cdots+n_s)^r}
\]
is convergent and its value is in the algebra of multiple zeta values.
\end{Cor}
Corollary \ref{zho} is actually the main result of \cite{zh} (From Theorem \ref{main} and Corollary \ref{rea}, we know that all the multiple zeta values in Corollary \ref{zho} are of weight $k_1+k_2+\cdots+k_s+r$ and depth $s$).  Yamamoto \cite{yama}
also treated this kind of series.

To calculate  iterated integrals on products of one variable multiple polylogarithms in general cases, we define the multiple Beta function.
\begin{Def}\label{mb}
For $N\geq 2$ and $\alpha_1,\alpha_2,\cdots,\alpha_N,\beta_1,\beta_2,\cdots,\beta_N>0$, define the multiple Beta function $B(\alpha_1,\beta_1; \alpha_2,\beta_2;\cdots, \alpha_N,\beta_N)$ as 
\[
\begin{split}
&\;\;\;\;\;\;B(\alpha_1,\beta_1;\alpha_2,\beta_2;\cdots; \alpha_N,\beta_N)\\
&=\mathop{\int\cdots\int}_{0<t_1<t_2<\cdots<t_N<1}t_1^{\alpha_1-1}(1-t_1)^{\beta_1-1}t_2^{\alpha_2-1}(1-t_2)^{\beta_2-1}\cdots t_N^{\alpha_N-1}(1-t_N)^{\beta_N-1}dt_1dt_2\cdots dt_N
\end{split}
\]
\end{Def}

Classical multiple zeta values can be viewed as evaluations of the  multiple Beta functions at special (limit) points. In fact, from the iterated integral representations of multiple zeta values we have 
\begin{prop}
For multiple zeta value $\zeta(k_1,k_2,\cdots,k_r)$, if we write $N=k_1+k_2+\cdots+k_r$, then
\[
\zeta(k_1,k_2,\cdots,k_r)=\mathop{{\mathrm{lim}}}_{\substack{(\alpha_i,\beta_i)\rightarrow (1,0)^+,\;i\in\{1,k_1+1,\cdots,k_1+\cdots+k_{r-1}+1\}   \\ (\alpha_i,\beta_i)\rightarrow (0,1)^+,\;i\notin\{1,k_1+1,\cdots,k_1+\cdots+k_{r-1}+1\}  }}B(\alpha_1,\beta; \alpha_2,\beta_2;\cdots, \alpha_N,\beta_N).
\]
\end{prop}

\begin{Thm}\label{beta} For $\alpha_1,\alpha_2,\cdots,\alpha_N,\beta_1,\beta_2,\cdots,\beta_N>0$,\\
$(i)$  For $1\leq i\leq N$, we have
\[
\begin{split}
&\;\;\;\;\;B(\alpha_1,\beta_1; \alpha_2,\beta_2;\cdots; \alpha_N,\beta_N)\\
&=B(\alpha_1,\beta_1;\cdots;\alpha_i+1,\beta_i;\cdots; \alpha_N,\beta_N)+B(\alpha_1,\beta_1;\cdots;\alpha_i,\beta_i+1;\cdots; \alpha_N,\beta_N);\\
\end{split}
\]
$(ii)$ $B(\alpha_1,\beta_1; \alpha_2,\beta_2;\cdots; \alpha_N,\beta_N)=B(\beta_N,\alpha_N;\beta_{N-1},\alpha_{N-1};\cdots;\beta_1,\alpha_1) $;\\
$(iii)$ 
\[
\begin{split}
&\;\;\;\;\frac{1}{\alpha_1}B(\alpha_1+1,\beta_1;\alpha_2,\beta_2;\cdots;\alpha_N,\beta_N)-\frac{1}{\beta_1}B(\alpha_1,\beta_1+1;\alpha_2,\beta_2;\cdots;\alpha_N,\beta_N)\\
&=-\frac{1}{\alpha_1\beta_1}B(\alpha_1+\alpha_2,\beta_1+\beta_2;\alpha_3,\beta_3;\cdots;\alpha_N,\beta_N),\\
\end{split}
\]
\[
\begin{split}
&\;\;\;\;\;\frac{1}{\alpha_i}B(\alpha_1,\beta_1;\cdots;\alpha_i+1,\beta_i;\cdots;\alpha_N,\beta_N)-\frac{1}{\beta_i}B(\alpha_1,\beta_1;\cdots;\alpha_i,\beta_i+1;\cdots;\alpha_N,\beta_N)\\
&=-\frac{1}{\alpha_i\beta_i}\bigg{( }B(\alpha_1,\beta_1;\cdots;\alpha_{i-1},\beta_{i-1};\alpha_i+\alpha_{i+1},\beta_i+\beta_{i+1};\alpha_{i+2},\beta_{i+2};\cdots; \alpha_N,\beta_N)\\
&\;\;\;\;\;\;\;\;\;\;\;\;\;\;\;\;\;\;-B(\alpha_1,\beta_1;\cdots;\alpha_{i-2},\beta_{i-2};\alpha_{i-1}+\alpha_i,\beta_{i-1}+\beta_i;\alpha_{i+1},\beta_{i+1};\cdots; \alpha_N,\beta_N)\bigg{)},\\
\end{split}
\]
\[
\begin{split}
&\;\;\;\;\;\frac{1}{\alpha_N}B(\alpha_1,\beta_1;\cdots;\alpha_{N-1},\beta_{N-1};\alpha_N+1,\beta_N)-\frac{1}{\beta_N}B(\alpha_1,\beta_1;\cdots; \alpha_{N-1},\beta_{N-1}; \alpha_{N},\beta_N+1)\\
&=\frac{1}{\alpha_N\beta_N}B(\alpha_1,\beta_1;\cdots; \alpha_{N-2},\beta_{N-2};\alpha_{N-1}+\alpha_N,\beta_{N-1}+\beta_N),\\
\end{split}
\]
$(iv)$ $B(\alpha_1,1;\alpha_2,1;\cdots;\alpha_N,1)=\frac{1}{\alpha_1(\alpha_1+\alpha_2)\cdots(\alpha_1+\alpha_2+\cdots+\alpha_N)}.$
\end{Thm}
\noindent{\bf Proof:}
(i) follows immediately from 
\[
\begin{split}
&\;\;\;\;\;t_1^{\alpha_1-1}(1-t_1)^{\beta_1-1}t_2^{\alpha_2-1}(1-t_2)^{\beta_2-1}\cdots t_N^{\alpha_N-1}(1-t_N)^{\beta_N-1}\\
&= t_1^{\alpha_1-1}(1-t_1)^{\beta_1-1}\cdots [t_i+(1-t_i)]t_{i}^{\alpha_{i}-1}(1-t_{i})^{\beta_{i}-1}\cdots t_N^{\alpha_N-1}(1-t_N)^{\beta_N-1}  \\
&=t_1^{\alpha_1-1}(1-t_1)^{\beta_1-1}\cdots t_i^{\alpha_i}(1-t_i)^{\beta_i-1}\cdots t_N^{\alpha_N-1}(1-t_N)^{\beta_N-1}  \\
&\;\;\;\;\;\;\;\;\;\;\;\;\;\;\;\;\;\;\;\;\;\;\;+t_1^{\alpha_1-1}(1-t_1)^{\beta_1-1}\cdots t_{i-1}^{\alpha_i}(1-t_i)^{\beta_i}\cdots t_N^{\alpha_N-1}(1-t_N)^{\beta_N-1}.
\end{split}
\]
Let $t_1=1-s_1,t_2=1-s_2,\cdots,t_N=1-s_N$, we have 
\[
\begin{split}
&\;\;\;\;\;\;B(\alpha_1,\beta_1;\alpha_2,\beta_2;\cdots; \alpha_N,\beta_N)\\
&=\mathop{\int\cdots\int}_{1>s_1>s_2>\cdots>s_N>0}(1-s_1)^{\alpha_1-1}s_1^{\beta_1-1}(1-s_2)^{\alpha_2-1}s_2^{\beta_2-1}\cdots (1-s_N)^{\alpha_N-1}s_N^{\beta_N-1}ds_1ds_2\cdots ds_N\\
&=\mathop{\int\cdots\int}_{0<s_1<s_2<\cdots<s_N<1}s_N^{\beta_N-1}(1-s_N)^{\alpha_N-1}\cdots s_2^{\beta_2-1}(1-s_2)^{\alpha_2-1}s_1^{\beta_1-1}(1-s_1)^{\alpha_1-1}ds_N\cdots ds_{2}ds_1\\
&=B(\beta_N,\alpha_N;\cdots;\beta_2,\alpha_2;\beta_1,\alpha_1).
\end{split}
\]
(iii) actually follows from the following simple formulas
\[
\frac{d}{dt_i}\left(\frac{1}{\alpha_i\beta_i}t_i^{\alpha_i}(1-t_i)^{\beta_i}  \right)=-\frac{1}{\alpha_i}t_i^{\alpha_i}(1-t_i)^{\beta_i-1}+\frac{1}{\beta_i}t_i^{\alpha_i-1}(1-t_i)^{\beta_i},1\leq i\leq N.
\]
(iv) follows from 
\[
\begin{split}
&\;\;\;\;\;B(\alpha_1,1;\alpha_2,1;\cdots;\alpha_N,1)\\
&=\mathop{\int\cdots\int}_{0<t_1<t_2<\cdots<t_N<1}t_1^{\alpha_1-1}t_2^{\alpha_2-1}\cdots t_N^{\alpha_N-1}dt_1dt_2\cdots dt_N\\
&=\frac{1}{\alpha_1}\mathop{\int\cdots\int}_{0<t_2<\cdots<t_N<1}t_2^{\alpha_1+\alpha_2-1}\cdots t_N^{\alpha_N-1}dt_2\cdots dt_N\\
&=\cdots\\
&=\frac{1}{\alpha_1}\cdot\frac{1}{\alpha_1+\alpha_2}\cdot\;\;\;\cdots \;\;\;\cdot\frac{1}{\alpha_1+\alpha_2+\cdots+\alpha_N}
\end{split}
\]
$\hfill\Box$\\

For $f_1(z),f_2(z),\cdots,f_N(z)\in {\bf{P}}_{\mathbb{Q}}^{\,\mathrm{log}}$, denote by 
\[
f_i(z)\,\omega_i(z)=\sum_{m,n\geq 1}a^i_{m,n}z^{m-1}(1-z)^{n-1}dz,1\leq i\leq N,
\]
then 
\[
\begin{split}
&\;\;\;\;\;\mathop{\int\cdots\int}_{0<z_1<z_2<\cdots<z_N<1}f_1(z_1)\omega_1(z_1)f_2(z_2)\omega_2(z_2)\cdots f_N(z_N)\omega_N(z_N)\\
&=\sum_{m_i,n_i\geq 1, 1\leq i\leq N}a^1_{m_1,n_1}a^2_{m_2,n_2}\cdots a^N_{m_N,n_N}\\
&\cdot\mathop{\int\cdots\int}_{0<z_1<z_2<\cdots<z_N<1}z_1^{m_1-1}(1-z_1)^{n_1-1}z_2^{m_2-1}(1-z_2)^{n_2-1}\cdots z_N^{m_N-1}(1-z)^{n_N-1}dz_1dz_2\cdots dz_N   \\
&= \sum_{m_i,n_i\geq 1, 1\leq i\leq N}a^1_{m_1,n_1}a^2_{m_2,n_2}\cdots a^N_{m_N,n_N} B(m_1,n_1;m_2,n_2;\cdots;m_N,n_N). \\
\end{split}
\]

The number $B(m_1,n_1;m_2,n_2;\cdots;m_N,n_N)$ can be calculated inductively by Theorem \ref{beta}.
From Theorem \ref{main}, the above integral is a $\mathbb{Q}$-linear combination of multiple zeta values if it is convergent. 
Unfortunately we don't know whether there exists a nice formula for the number $$B(m_1,n_1;m_2,n_2;\cdots;m_N,n_N)$$ so far.

\section*{Acknowledgements}
The author wants to thank Professor Qingchun Tian for his careful reading and useful comments to improve this paper. The author also thanks Professor Henrik Bachmann for the notification of Seki  and Yamamoto's paper \cite{SY}. The author would like to express his sincere gratitude to the anonymous referee for  his/her detailed comments to improve this paper.

\end{document}